\documentclass[12pt,letterpaper,oneside]{amsart}
\usepackage{amsmath,amssymb}
\usepackage[left=1in,textwidth=6.3in,top=1.27in,textheight=8.4in,footskip=33pt]{geometry}
\usepackage{tikz}
\usepackage[colorlinks=true,citecolor=blue,linkcolor=blue,urlcolor=black,pagebackref]{hyperref}

\hypersetup{pdftitle={Volume difference inequalities for central sections of origin-symmetric convex bodies}}
\numberwithin{equation}{section}
\newtheorem{theorem}{Theorem}[section]
\newtheorem{proposition}[theorem]{Proposition}
\newtheorem{lemma}[theorem]{Lemma}
\theoremstyle{definition}
\newtheorem{gkquestion}{Question}[section]

\DeclareMathOperator{\conv}{conv}
\newcommand{\R}{\mathbb{R}}
\newcommand{\B}{B_2}
\newcommand{\sph}{\mathbb{S}^{n-1}}
\newcommand{\SF}{S_F}
\newcommand{\dd}{\,d}
\newcommand{\E}{\mathbb E}
\newcommand{\T}{\mathsf T}

\DeclareMathOperator{\supp}{supp}
\DeclareMathOperator{\Lip}{Lip}

\newcommand{\nnfootnote}[1]{%
  \begingroup
  \renewcommand{\thefootnote}{}%
  \footnotetext[0]{#1}%
  \endgroup
}
\newcommand*{\manuscriptaffiliation}{}
\makeatletter
\renewcommand{\@settitle}{%
  \global\topskip\normaltopskip
  {\centering\large\bfseries\@title\par}%
}
\renewcommand{\@setauthors}{%
  \begingroup
  \author@andify\authors
  \centering\medskip
  {\normalsize\bfseries\authors\par}%
  \ifx\manuscriptaffiliation\@empty\else
    \medskip
    {\small\manuscriptaffiliation\par}%
  \fi
  \endgroup
}
\renewcommand{\@setabstracta}{%
  \ifvoid\abstractbox\else
    \vskip8pt
    \box\abstractbox
    \prevdepth\z@
  \fi
}
\renewenvironment{abstract}{%
  \global\setbox\abstractbox=\vtop\bgroup
  \normalfont\normalsize
  \list{}{%
    \leftmargin28pt \rightmargin\leftmargin
    \labelwidth\z@ \labelsep\z@ \itemindent\z@ \listparindent\z@
    \topsep\z@ \partopsep\z@ \parsep\z@
  }%
  \item[]\textbf{Abstract:}\enspace\ignorespaces
}{%
  \endlist\egroup
}
\renewcommand{\section}{%
  \@startsection{section}{1}{\z@}%
  {\linespacing\@plus.4\linespacing\@minus.2\linespacing}%
  {.7\linespacing}{\normalfont\centering\bfseries}%
}
\renewcommand{\subsection}{%
  \@startsection{subsection}{2}{\z@}%
  {.8\linespacing\@plus.2\linespacing}%
  {.25\linespacing}{\normalfont\bfseries}%
}
\renewcommand{\@bibtitlestyle}{\specialsection*{\scshape\refname}}
\makeatother

\title{}

\begin{document}
\begin{center}
{\Large \bf  Volume difference inequalities for central sections\\
of origin-symmetric convex bodies}
\end{center}
\vskip 10pt
\begin{center}
{ Jiangyan\ Tao}\\~~ \\
\small{School of Mathematical Sciences, Zhejiang Normal University, Zhejiang, 321004, China}
\end{center}
\vskip 5pt
\nnfootnote{E-mail address: taojiangyan@zjnu.edu.cn.}

\begin{abstract}
Giannopoulos and Koldobsky~\cite{GK} asked whether, for origin-symmetric
convex bodies $L\subset K\subset\R^n$, the difference
$|K|^{(n-k)/n}-|L|^{(n-k)/n}$ is bounded by $C^k$ times the largest
difference of their central section volumes in codimension $k$, with a
single absolute constant $C$. We prove that the optimal constant
$\gamma_{n,k}$ for fixed $n,k$ satisfies $\gamma_{n,k}\ge c\sqrt{n/k}$
for an absolute $c>0$. For $k=1$, this lower bound grows as $\sqrt n$
with the dimension, giving a negative answer to the question.

\smallskip
\noindent\textbf{2020 Mathematics Subject Classification:} 52A20, 52A40.

\smallskip
\noindent\textbf{Keywords:} Volume difference, central sections, slicing inequality.
\end{abstract}

\maketitle
\vspace{-8pt}

\section{Introduction}

Let \(K,L\subset\R^n\) be origin-symmetric convex bodies, that is, compact
convex sets with nonempty interior satisfying \(K=-K\) and \(L=-L\).
For \(1\le k<n\), let \(G_{n,n-k}\) denote the set of all
\((n-k)\)-dimensional linear subspaces of \(\R^n\). For each
\(F\in G_{n,n-k}\), the intersection \(K\cap F\) is a central section
of \(K\) of codimension \(k\). We write \(|K|\) for the \(n\)-dimensional
volume of \(K\) and \(|K\cap F|\) for the \((n-k)\)-dimensional volume
of \(K\cap F\) measured in \(F\). We denote the Euclidean unit ball in
\(\R^d\) by \(B_2^d\). The letters \(c\) and \(C\) denote positive
absolute constants whose values may vary.

The classical Busemann--Petty problem~\cite{BP} asks whether the
inequalities
\[
|K\cap H|\le |L\cap H|
\qquad\text{for every }H\in G_{n,n-1}
\]
imply the volume inequality
\[
|K|\le |L|.
\]
The answer is affirmative for \(n\le4\) and negative for
\(n\ge5\)~\cite{GKS}. Zhang~\cite{Zhang} posed the corresponding
question for central sections of higher codimension. For further work
on these problems, see~\cite{GardnerBP3,KoldobskyBP,LarmanRogersBP,LLX,LutwakIntersectionBP,ZhangBP4}.

The stability question allows the volume of each section of \(K\) to
exceed that of the corresponding section of \(L\) by at most
\(\varepsilon\). For fixed \(n\) and \(k\), one asks whether there is a
constant \(A=A(n,k)\) such that, for all \(K,L\) and all
\(\varepsilon\ge 0\), the section bound
\[
\max_{F\in G_{n,n-k}}
\bigl(|K\cap F|-|L\cap F|\bigr)\le \varepsilon
\]
implies the volume bound
\[
|K|^{(n-k)/n}-|L|^{(n-k)/n}\le A\varepsilon.
\]
When \(k=1\) and \(\varepsilon=0\), the volume bound reduces to
\(|K|\le |L|\), so one recovers the classical Busemann--Petty problem.
Its counterexamples in dimensions \(n\ge 5\) therefore rule out this
stability estimate for arbitrary pairs of origin-symmetric convex bodies
in the hyperplane case.

For every \(F\in G_{n,n-k}\),
\[
 |K\cap F|-|L\cap F|
 = |(K\setminus L)\cap F|-|(L\setminus K)\cap F|.
\]
For a general pair, the two terms on the right can both be large even
when their difference is small. If \(L\subset K\), the second term
vanishes, and the smallest value of \(\varepsilon\) in the section
bound is
\[
\varepsilon
=\max_{F\in G_{n,n-k}} |(K\setminus L)\cap F|.
\]
Substituting this value of \(\varepsilon\) into the stability estimate
and writing \(A=C^k\) leads to the following question, posed by
Giannopoulos and Koldobsky~\cite[Question~1.1]{GK}.

\begin{gkquestion}\label{ques:volume-difference}
Does there exist an absolute constant \(C>0\) such that
\begin{equation}\label{eq:volume-difference}
|K|^{(n-k)/n}-|L|^{(n-k)/n}
\le C^k\max_{F\in G_{n,n-k}}
\bigl(|K\cap F|-|L\cap F|\bigr)
\end{equation}
for all integers \(n\ge 2\) and \(1\le k<n\) and all origin-symmetric
convex bodies \(L\subset K\subset\R^n\)?
\end{gkquestion}

For each fixed \(n,k\), let \(\gamma_{n,k}\) be the least constant \(C\) in
\eqref{eq:volume-difference}. Gardner~\cite[report on Problem~8.3]{GardnerUpdate}
discusses this question under the name \emph{volume difference conjecture}.
The exponent \((n-k)/n\) is dictated by scaling. 
For \(L=tK\) with \(0<t<1\), the identity \(tK\cap F=t(K\cap F)\) allows us
to cancel \(1-t^{n-k}\) on both sides of \eqref{eq:volume-difference}. This gives
the slicing inequality for a single body,
\[
 |K|^{(n-k)/n}\le C^k\max_{F\in G_{n,n-k}}|K\cap F|.
\]
For \(k=1\), the existence of an absolute \(C\) in this last inequality is
the origin-symmetric case of Bourgain's slicing problem~\cite{Bourgain}.
Klartag and Lehec~\cite[Theorem~1.1]{KLslicing} proved that every convex body
has an affine hyperplane section of volume at least \(c|K|^{(n-1)/n}\).
For an origin-symmetric body, the central section has the largest volume
among all hyperplane sections parallel to it. Indeed, by Brunn's concavity
theorem the \((n-1)\)-st root of the section volume is concave in the signed
distance from the origin on the interval of nonempty sections. Symmetry makes
this function even, so it is maximal at zero~\cite[Section~1.3]{BGVV}.

The first estimates for \eqref{eq:volume-difference} used intersection bodies.
Koldobsky~\cite[Theorem~1]{KoldobskyStability} proved the inequality with
\(C=1\) when \(k=1\) and \(K\) is an intersection body. Koldobsky and
Ma~\cite[Theorem~1]{KM} obtained the coefficient
\(C^k=|B_2^n|^{(n-k)/n}/|B_2^{n-k}|\) when \(K\) is a generalized
\(k\)-intersection body. See also the survey of Giannopoulos, Koldobsky, and
Zvavitch~\cite[Section~7]{GKZSurvey}. Giannopoulos and
Koldobsky~\cite[Theorem~1.2]{GK} extended this estimate to arbitrary
origin-symmetric convex bodies. Their bound multiplies the preceding
coefficient by the \(k\)-th power of the outer volume ratio distance
to the class of generalized \(k\)-intersection bodies. This distance is the
infimum of \((|D|/|K|)^{1/n}\) over bodies \(D\) in that class containing
\(K\). It is uniformly bounded when \(K\) is unconditional, that is,
invariant under changes of coordinate signs.

For general origin-symmetric convex bodies the known upper bounds give
\[
 \gamma_{n,k}\le\sqrt n,\qquad
 \gamma_{n,k}\le C\min\left\{\sqrt{n-k},\,
    \sqrt{\frac nk}\left(\log\frac{en}{k}\right)^{3/2}\right\}.
\]
The bound \(\sqrt n\) follows from the outer volume ratio estimate given
by John's ellipsoid theorem~\cite[Theorem~1.2 and the subsequent discussion]{GK}.
For the logarithmic bound, Giannopoulos and
Koldobsky~\cite[Corollary~1.3]{GK} use the outer volume ratio estimate of
Koldobsky, Paouris, and Zymonopoulou~\cite{KPZ}, with the logarithmic exponent
from~\cite[Proposition~1]{KoldobskySlicing}. The bound \(C\sqrt{n-k}\)
follows by applying the section inequality of Chasapis, Giannopoulos, and
Liakopoulos~\cite[Theorem~1.1]{CGL} to the indicator of \(K\setminus L\),
as in~\cite[Remark~2.3, equation~(2.23)]{GK}, and using
\[
 |K|^{(n-k)/n}-|L|^{(n-k)/n}
 \le \frac{|K|-|L|}{|K|^{k/n}}.
\]

The known upper bounds leave open the possibility that a single absolute
constant suffices in \eqref{eq:volume-difference} for all \(n\) and \(k\).
To rule this out, it suffices to find a lower bound for \(\gamma_{n,1}\)
that tends to infinity with \(n\). We shall prove more generally that
\(\gamma_{n,k}\) tends to infinity whenever \(k/n\to 0\).

Our first result is a new characterization of the optimal constant
\(\gamma_{n,k}\) in the volume difference inequality
\eqref{eq:volume-difference}. In Proposition~\ref{prop:density} we prove that
\[
\gamma_{n,k}^{\,k}
=\frac{n-k}{n}\sup_{K,f}
\frac{1}{\displaystyle
	\max_{F\in G_{n,n-k}}\int_{K\cap F} f(x)\,dx_F},
\]
where \(K\) ranges over origin-symmetric convex bodies of volume one in
\(\R^n\), \(f\) ranges over even continuous probability densities on \(K\),
and \(dx_F\) denotes Lebesgue measure on \(F\). Thus, to obtain a lower
bound for the volume difference constant, it suffices to find one pair
\(K,f\) for which all central section integrals are small.

Such examples arose in the study of slicing inequalities for
measures~\cite{KoldobskySlicing}. For \(k=1\) and \(n\ge3\), Klartag and
Koldobsky~\cite[Theorem~1.1]{KKexample} constructed a density for which
every hyperplane section integral is at most \(C\sqrt{(\log\log n)/n}\).
Klartag and Livshyts~\cite[Theorem~1.1 in the corrected version]{KLrounding}
improved this bound to \(C/\sqrt n\). Applying our characterization to
the first example already gives \(\gamma_{n,1}\to\infty\) as \(n\to\infty\),
which answers Question~\ref{ques:volume-difference} in the negative. Applying it to the second
example gives the sharper bound \(\gamma_{n,1}\ge c\sqrt n\), matching
the known upper bound up to an absolute constant.

For general codimension \(k\), Haddad~\cite[Theorem~1.4]{Haddad}
constructed a density for which every section integral is at most
\((C\sqrt{k\log n/n})^k\). Substituting this estimate into
Proposition~\ref{prop:density} gives
\[
\gamma_{n,k}\ge
c\left(\frac{n-k}{n}\right)^{1/k}\sqrt{\frac{n}{k\log n}} .
\]
For \(2\le k\le n/2\), we remove the logarithmic loss in this bound
by constructing a density with smaller section integrals in
Lemma~\ref{lem:small-sections}. For \(k>n/2\), two concentric balls
suffice. Together these arguments give the following theorem.

\begin{theorem}\label{thm:main}
There is an absolute \(c>0\) such that, for all \(n\ge2\) and \(1\le k<n\),
\begin{equation}\label{eq:bounds}
 \gamma_{n,k}\ge c\sqrt{\frac nk}.
\end{equation}
\end{theorem}

For hyperplanes, Theorem~\ref{thm:main} and the upper bound above give
\(c\sqrt n\le\gamma_{n,1}\le\sqrt n\). Thus any proposed absolute constant
in Question~\ref{ques:volume-difference} fails in some dimension already for \(k=1\). Geometrically,
the ratio of the loss measured by \(|K|^{(n-1)/n}-|L|^{(n-1)/n}\) to the
largest loss in a central hyperplane section can grow without bound. The
single-body slicing theorem controls this ratio for the special pairs \(L=tK\).
More generally, for any sequence \((n_j,k_j)\) with \(1\le k_j<n_j\),
the constants \(\gamma_{n_j,k_j}\) are uniformly bounded if and only if
the ratios \(k_j/n_j\) are bounded away from zero.
Indeed, if \(\gamma_{n_j,k_j}\le M\) for every \(j\), then
\eqref{eq:bounds} gives \(k_j/n_j\ge (c/M)^2\) for every \(j\).
Conversely, if \(k_j/n_j\ge\alpha\) for some \(\alpha>0\) and every \(j\),
then the logarithmic upper bound gives
\[
\gamma_{n_j,k_j}
\le C\alpha^{-1/2}\bigl(\log(e/\alpha)\bigr)^{3/2},
\]
which is independent of \(j\).

The first step is to relate section integrals of densities to differences
of section volumes for convex bodies \(L\subset K\). By perturbing a polar
body, we can control how the radial lengths change while preserving
convexity. Together with an approximation by densities concentrated near
the boundary, this gives the formula for \(\gamma_{n,k}^{\,k}\) in
Proposition~\ref{prop:density}.

For \(k\le n/2\), we use convolution at several scales to construct a
probability density whose integral is small on every central section of
codimension \(k\). Maurey's covering lemma gives the required bound on
the volume of the supporting convex body. For \(k>n/2\), the lower bound
follows by taking \(K\) and \(L\) to be concentric Euclidean balls.

The paper is organized as follows. Section~\ref{sec:preliminaries} collects the polar integration formulas,
estimates for the projections of Euclidean balls, and the covering lemma for
finite convex hulls. In Section~\ref{sec:proof} we prove the density identity
of Proposition~\ref{prop:density}, construct a density with small central
section integrals in Lemma~\ref{lem:small-sections}, and conclude with the
proof of Theorem~\ref{thm:main}.

\vskip 20pt
\section{Preliminaries}\label{sec:preliminaries}

Throughout, \(n\ge2\) and \(1\le k<n\). We write \(\omega_d=|B_2^d|\)
and denote the Euclidean norm by \(|\cdot|\). Let
\(\sph=\{\theta\in\R^n:|\theta|=1\}\) and set \(\SF=\sph\cap F\)
for \(F\in G_{n,n-k}\).

For a convex body \(K\) with \(0\in\operatorname{int}K\), the radial
function is defined by
\[
\rho_K(\theta)=\max\{r\ge 0 : r\theta\in K\},
\qquad \theta\in\sph.
\]
It is positive and continuous on \(\sph\).

For a convex body \(K\) containing the origin in its interior, write
\[
 K^\circ=\{y\in\R^n:\langle x,y\rangle\le1
              \ \text{for every }x\in K\},\qquad
 h_K(v)=\max_{x\in K}\langle x,v\rangle .
\]
For \(v\in\sph\) and \(r\ge0\), the point \(rv\) belongs to \(K^\circ\)
if and only if \(r h_K(v)\le1\). Thus \(\rho_{K^\circ}(v)=1/h_K(v)\).
We also have \((K^\circ)^\circ=K\). Indeed, the inclusion
\(K\subseteq(K^\circ)^\circ\) follows from the definition.
If \(x\notin K\), the separation theorem gives \(v\in\sph\) such that
\(\langle x,v\rangle>h_K(v)>0\). Then \(v/h_K(v)\in K^\circ\) and
\(\langle x,v/h_K(v)\rangle>1\), so \(x\notin(K^\circ)^\circ\).

We will use two further consequences of the definitions. Inclusion
reverses under polarity, and for the Minkowski sum
\(K_1+tK_2=\{x+ty:x\in K_1,\,y\in K_2\}\), one has
\(h_{K_1+tK_2}=h_{K_1}+t h_{K_2}\) for \(t\ge0\).

Let \(\sigma\) and \(\sigma_F\) denote the surface area measures on
\(\sph\) and \(\SF\), respectively, with total masses \(n\omega_n\)
and \((n-k)\omega_{n-k}\). Polar coordinates give
\begin{equation}\label{eq:polar}
|K|=\frac1n\int_{\sph}\rho_K(\theta)^n\dd\sigma(\theta),
\qquad
|K\cap F|=\frac1{n-k}\int_{\SF}
\rho_K(\theta)^{n-k}\dd\sigma_F(\theta).
\end{equation}
More generally, for every continuous real-valued function \(f\) on \(K\),
\begin{align}
\int_K f(x)\dd x
&=\int_{\sph}\int_0^{\rho_K(\theta)}
f(r\theta)\,r^{n-1}\dd r\dd\sigma(\theta),
\label{eq:polar-density}\\
\int_{K\cap F} f(x)\dd x_F
&=\int_{\SF}\int_0^{\rho_K(\theta)}
f(r\theta)\,r^{n-k-1}\dd r\dd\sigma_F(\theta),
\label{eq:section-density}
\end{align}
where \(\dd x_F\) denotes Lebesgue measure on \(F\).

For continuous \(h\) on \(\sph\), set
\[
M(h)=\max_{F\in G_{n,n-k}}\int_{\SF} h\dd\sigma_F .
\]
The maximum is attained because the section integral depends continuously
on \(F\) and \(G_{n,n-k}\) is compact. If \(h\ge0\) and \(h\not\equiv0\),
choose \(F\) containing a point where \(h>0\). By continuity, the integral
over \(\SF\) is positive, so \(M(h)>0\).

To estimate the density of the projection of the uniform measure on a
Euclidean ball, we use the volume bounds
\begin{equation}\label{eq:ball-volume}
\frac2{\sqrt d}\le\omega_d^{1/d}\le\frac{\sqrt{2\pi e}}{\sqrt d},
\qquad d\ge1.
\end{equation}
The cube \([-1/\sqrt d,1/\sqrt d]^d\) is
contained in \(\B^d\), giving the lower bound. For the upper bound, integrate
\(e^{-|x|^2/2}\) on \(\sqrt d\,\B^d\), where it is at least \(e^{-d/2}\).
Comparison with its integral on \(\R^d\) gives
\(e^{-d/2}d^{d/2}\omega_d\le(2\pi)^{d/2}\).

Let \(1\le k<n\) be integers and let \(r>0\). Write \(\mu_r\) for the
uniform probability measure on the Euclidean ball \(rB_2^n\).
Let \(U\) be a \(k\times n\) matrix with orthonormal rows. The projection
of \(\mu_r\) under \(U\) has density \(q_r\) on \(\mathbb{R}^k\).
Thus, for every bounded Borel function \(h\) on \(\mathbb{R}^k\),
\[
\int_{rB_2^n} h(Ux)\,d\mu_r(x)
=\int_{\mathbb{R}^k} h(z)\,q_r(z)\,dz.
\]

In the orthogonal coordinates \(x=U^{\mathsf T}z+y\), with \(y\in\ker U\),
we have \(Ux=z\) and \(|x|^2=|z|^2+|y|^2\). For \(|z|<r\), the section
of \(rB_2^n\) at \(Ux=z\) is therefore an \((n-k)\)-dimensional ball
of radius \(\sqrt{r^2-|z|^2}\). For \(|z|\ge r\), its
\((n-k)\)-dimensional volume is zero. Since \(dx=dz\,dy\), Fubini's
theorem gives \(q_r(z)\) by dividing the section volume by
\(|rB_2^n|=\omega_n r^n\). Hence
\begin{equation}\label{eq:marginal}
q_r(z)=\frac{\omega_{n-k}}{r^k\omega_n}
\Big(1-\frac{|z|^2}{r^2}\Big)_+^{(n-k)/2},
\qquad z\in\mathbb{R}^k,
\end{equation}
where \(\omega_d=|B_2^d|\) and \(u_+=\max\{u,0\}\).
This formula gives the same density for every such \(U\). For \(r=1\),
it follows by taking \(\beta=0\) in the projection formula for beta
distributions given by Kabluchko, Temesvari, and
Th\"ale~\cite[Lemma~4.3(a)]{KTTBeta}.

\begin{lemma}\label{lem:ball-marginal}
For \(r>0\) and \(1\le k<n\), the density \(q_r\) in~\eqref{eq:marginal}
satisfies
\[
\|q_r\|_\infty=q_r(0)\le\left(\frac{\sqrt n}{r}\right)^k.
\]
If \(n-k\ge2\), its least Lipschitz constant on \(\mathbb{R}^k\) satisfies
\[
\operatorname{Lip}(q_r)\le\frac{\sqrt{2n}}{r}\,\|q_r\|_\infty.
\]
\end{lemma}

\begin{proof}
By~\eqref{eq:marginal}, the density attains its maximum at the origin.
To estimate \(q_r(0)\),
we use \(\omega_{j-1}/\omega_j\le\sqrt j\) for \(j\ge2\).
Alzer~\cite{AlzerBall} gives sharper bounds for this volume ratio.
The bound needed here follows by slicing \(B_2^j\) perpendicular to
its last coordinate. For \(|s|\le1/\sqrt{2j}\), Bernoulli's inequality
gives \((1-s^2)^{j-1}\ge 1-(j-1)s^2\ge 1/2\). Hence
\[
\frac{\omega_j}{\omega_{j-1}}
=\int_{-1}^{1}(1-s^2)^{(j-1)/2}\,ds
\ge \frac{2}{\sqrt{2j}}\cdot\frac1{\sqrt2}
=\frac1{\sqrt j}.
\]
Multiplying the reciprocal bounds for \(j=n-k+1,\ldots,n\) gives
\[
q_r(0)=\frac{1}{r^k}
\prod_{j=n-k+1}^{n}\frac{\omega_{j-1}}{\omega_j}
\le \frac{1}{r^k}\prod_{j=n-k+1}^{n}\sqrt j
\le \left(\frac{\sqrt n}{r}\right)^k.
\]
	
For the Lipschitz estimate, assume \(n-k\ge2\) and set
\(g(t)=(1-t^2)_+^{(n-k)/2}\) for \(t\ge0\). For \(0<t<1\),
\[
|g'(t)|^2=(n-k)^2 t^2(1-t^2)^{n-k-2}
\le\frac{(n-k)^2}{n-k-1}\le 2(n-k).
\]
Indeed, \((n-k-1)t^2(1-t^2)^{n-k-2}\) is one of the nonnegative terms
in the binomial expansion of \(\bigl(t^2+(1-t^2)\bigr)^{n-k-1}=1\).
The mean value theorem gives a Lipschitz constant of \(\sqrt{2(n-k)}\)
for \(g\) on \([0,1]\). Since \(g(1)=0\) and \(g(t)=0\) for \(t\ge1\),
the same bound holds on \([0,\infty)\). Now
\(q_r(z)=q_r(0)\,g(|z|/r)\) and \(\bigl||z|-|w|\bigr|\le |z-w|\),
so \(\operatorname{Lip}(q_r)\le q_r(0)\sqrt{2(n-k)}/r\), as claimed.
When \(n-k=1\), the derivative of \(\sqrt{1-t^2}\) is unbounded near
\(t=1\), which explains the restriction \(n-k\ge2\).
\end{proof}

To estimate the volume of a finite convex hull, we use Maurey's
approximation lemma~\cite[Lemma~2]{Pisier} to cover it by balls centered
at averages of the given points. Counting multisets bounds the number of
centers. For a greedy approach to the approximation, see~\cite[Section~VIII]{Barron}.

\begin{lemma}\label{lem:convex-hull-covering}
For a finite nonempty set $V\subset rB_2^n$, write
$N=\operatorname{card}V$. For every integer $\ell\ge1$ there is
a finite set $D$ such that
\[
\conv V\subset D+\frac{r}{\sqrt\ell}B_2^n,\qquad
\operatorname{card}D\le \binom{N+\ell-1}{\ell}.
\]
\end{lemma}

\begin{proof}
Write \(V=\{v_1,\ldots,v_N\}\) and set
\[
D=\Big\{\frac1\ell\sum_{j=1}^{\ell}v_{i_j}:
1\le i_1\le\cdots\le i_\ell\le N\Big\}.
\]
Each index tuple is uniquely determined by the \(\ell\)-element set
\(\{i_j+j-1:1\le j\le\ell\}\subset\{1,\ldots,N+\ell-1\}\). Hence
\[
\operatorname{card}D\le\binom{N+\ell-1}{\ell}.
\]
	
Fix \(x\in\conv V\) and write
\[
x=\sum_{i=1}^N p_i v_i,\qquad p_i\ge0,\qquad \sum_{i=1}^N p_i=1.
\]
Let \(X_1,\ldots,X_\ell\) be independent random vectors with
\(\mathbb P(X_j=v_i)=p_i\). Then \(\E X_j=x\) and
\(\ell^{-1}\sum_{j=1}^{\ell}X_j\in D\) almost surely.
Since the vectors \(X_j-x\) are independent and centered,
\begin{align*}
\min_{y\in D}|x-y|^2
\le \E\Big|\frac1\ell\sum_{j=1}^{\ell}(X_j-x)\Big|^2
=\frac1{\ell^2}\sum_{j=1}^{\ell}\E|X_j-x|^2
=\frac1\ell\Big(\sum_{i=1}^N p_i|v_i|^2-|x|^2\Big)
\le \frac{r^2}{\ell},
\end{align*}
where the last inequality follows from \(|v_i|\le r\).
Thus \(\conv V\subset D+(r/\sqrt{\ell})B_2^n\).
\end{proof}

\vskip 20pt
\section{Proof of the main theorem}\label{sec:proof}

To obtain a lower bound for \(\gamma_{n,k}\), we need convex bodies
\(L\subset K\) for which the quotient of the volume difference by the
largest section difference is large. The following proposition expresses
the optimal constant in terms of integrals of a density on a single
convex body.

\begin{proposition}\label{prop:density}
Fix integers \(n\ge2\) and \(1\le k<n\). Let \(K\) range over
origin-symmetric convex bodies in \(\R^n\). For each \(K\), let \(f\)
range over continuous, even, nonnegative functions on \(K\) satisfying
\(\int_K f(x)\dd x>0\). Then
\begin{equation}\label{eq:density}
\gamma_{n,k}^{\,k}
=\frac{n-k}{n}\sup_{K,f}
\frac{\displaystyle\int_K f(x)\dd x}
{\displaystyle |K|^{k/n}
\max_{F\in G_{n,n-k}}\int_{K\cap F} f(x)\dd x}.
\end{equation}
\end{proposition}

\begin{proof}
\textbf{Step~1.} We use polar coordinates and H\"older's inequality
to bound \(\gamma_{n,k}^{\,k}\) from above by a supremum over functions
on the unit sphere.
	
For \(L\subsetneq K\), set \(h=\rho_K^{\,n-k}-\rho_L^{\,n-k}\).
Then \(h\) is continuous, even, and nonnegative. Since the inclusion
is strict, \(h(\theta)>0\) for some \(\theta\in\sph\).
By~\eqref{eq:polar}, for every \(F\in G_{n,n-k}\),
\[
|K\cap F|-|L\cap F|
=\frac{1}{n-k}\int_{\SF} h\,\dd\sigma_F .
\]
If \(F\) contains \(\theta\), then \(h\) is positive on a neighborhood
of \(\theta\) in \(\SF\), so this section difference is positive.
Since both sides of~\eqref{eq:volume-difference} vanish when \(L=K\),
the definition of \(\gamma_{n,k}\) gives
\[
\gamma_{n,k}^{\,k}
=\sup_{L\subsetneq K}
\frac{|K|^{(n-k)/n}-|L|^{(n-k)/n}}
{\displaystyle\max_{F\in G_{n,n-k}}
\bigl(|K\cap F|-|L\cap F|\bigr)} .
\]
	
H\"older's inequality with conjugate exponents \(n/k\) and \(n/(n-k)\) gives
\[
\int_{\sph}\rho_K^k\rho_L^{n-k}\dd\sigma
\le \Big(\int_{\sph}\rho_K^n\dd\sigma\Big)^{k/n}
\Big(\int_{\sph}\rho_L^n\dd\sigma\Big)^{(n-k)/n}
=n|K|^{k/n}|L|^{(n-k)/n}.
\]
Subtracting this estimate from \(\int_{\sph}\rho_K^n\dd\sigma=n|K|\)
and using \(h=\rho_K^{\,n-k}-\rho_L^{\,n-k}\), we obtain
\[
\int_{\sph}\rho_K^k h\dd\sigma
\ge n|K|^{k/n}\bigl(|K|^{(n-k)/n}-|L|^{(n-k)/n}\bigr).
\]
Divide the preceding inequality by \(n|K|^{k/n}\) and by the largest
section difference, which is positive for \(L\subsetneq K\).
Using~\eqref{eq:polar} to express this difference in terms of \(h\), we obtain
\begin{equation}\label{eq:angular-upper}
\frac{|K|^{(n-k)/n}-|L|^{(n-k)/n}}
{\displaystyle\max_{F\in G_{n,n-k}}
\bigl(|K\cap F|-|L\cap F|\bigr)}
\le \frac{n-k}{n}\,
\frac{\displaystyle\int_{\sph}\rho_K^k h\dd\sigma}
{\displaystyle |K|^{k/n}
\max_{F\in G_{n,n-k}}\int_{\SF}h\dd\sigma_F}.
\end{equation}
Taking the supremum over \(K\) and \(L\subsetneq K\) gives
\(\gamma_{n,k}^{\,k}\) on the left. For each \(K\), allowing \(h\)
on the right to range over all continuous, even, nonnegative functions
on \(\sph\) with \(h\not\equiv0\) can only increase the supremum.
	
\textbf{Step~2.} We show that for each \(K\) and \(h\) allowed in
Step~1, the right-hand side of~\eqref{eq:angular-upper} is a limit of
volume-difference quotients.
	
Fix an origin-symmetric convex body \(K\) and a continuous, even,
nonnegative function \(h\not\equiv0\) on \(\sph\). We seek convex
bodies \(L_t\subsetneq K_t\) such that
\[
\rho_{K_t}^{\,n-k}-\rho_{L_t}^{\,n-k}=t^2h+o(t^2)
\]
uniformly on \(\sph\) as \(t\downarrow0\). Since the derivative of
\(s^{-(n-k)}\) at \(s=1/\rho_K(v)\) is
\(-(n-k)\rho_K(v)^{n-k+1}\), we choose
\[
a(v)=\frac{h(v)}{(n-k)\rho_K(v)^{n-k+1}},
\qquad v\in\sph.
\]
For \(t>0\), let
\[
D_t=\conv\{(1+t a(u))u:u\in\sph\}.
\]
The function \(a\) is continuous, even, and nonnegative. Thus \(D_t\)
is compact and symmetric about the origin. For every \(u\in\sph\),
the segment joining \(\pm(1+t a(u))u\) contains \([-u,u]\).
Hence \(\B^n\subseteq D_t\).
	
We first determine the support function of \(D_t\) to first order in \(t\).
Fix \(v\in\sph\) and choose \(u\in\sph\) maximizing
\((1+t a(u))\langle u,v\rangle\). This maximum equals \(h_{D_t}(v)\),
since taking the convex hull preserves the maximum of a linear
functional. Comparing with the value at \(u=v\) gives
\[
1+t a(v)
\le h_{D_t}(v)
=(1+t a(u))\langle u,v\rangle
\le 1+t a(u).
\]
In particular, \(\langle u,v\rangle\ge(1+t\|a\|_\infty)^{-1}\), so
\[
|u-v|^2=2(1-\langle u,v\rangle)\le 2t\|a\|_\infty.
\]
It follows that
\[
0\le \frac{h_{D_t}(v)-1}{t}-a(v)\le |a(u)-a(v)|.
\]
Since \(a\) is uniformly continuous on \(\sph\) and the bound on
\(|u-v|\) is independent of \(v\), the right-hand side tends to zero
uniformly in \(v\). We have proved
\begin{equation}\label{eq:support-variation}
h_{D_t}(v)=1+t a(v)+o(t)
\qquad\text{uniformly for }v\in\sph.
\end{equation}
	
Now define
\[
K_t=(K^\circ+t\B^n)^\circ,
\qquad
L_t=(K^\circ+tD_t)^\circ.
\]
These are origin-symmetric convex bodies. Since \(\B^n\subseteq D_t\)
and polarity reverses inclusion, \(L_t\subseteq K_t\).
The identity \(\rho_{D^\circ}=1/h_D\) and the additivity of support
functions under Minkowski addition give
\[
\frac{1}{\rho_{K_t}(v)}=\frac{1}{\rho_K(v)}+t,
\qquad
\frac{1}{\rho_{L_t}(v)}=\frac{1}{\rho_K(v)}+t h_{D_t}(v).
\]
At any \(v\) with \(h(v)>0\), we have \(h_{D_t}(v)\ge1+t a(v)>1\).
Thus \(\rho_{L_t}(v)<\rho_{K_t}(v)\), proving \(L_t\subsetneq K_t\).
Moreover, by~\eqref{eq:support-variation},
\[
\frac{1}{\rho_{L_t}}-\frac{1}{\rho_{K_t}}
=t(h_{D_t}-1)=t^2a+o(t^2)
\]
uniformly on \(\sph\).
	
For \(d\in\{n-k,n\}\), we have
\[
\rho_{K_t}(v)^d-\rho_{L_t}(v)^d
=d\int_{1/\rho_{K_t}(v)}^{1/\rho_{L_t}(v)}s^{-d-1}\dd s.
\]
Both endpoints converge uniformly to \(1/\rho_K(v)\), which is bounded
away from zero on \(\sph\). Thus, throughout the interval of integration,
\(s^{-d-1}=\rho_K(v)^{d+1}+o(1)\) uniformly in \(v\). Multiplying by
the length of the interval and substituting the definition of \(a\)
gives
$$
\rho_{K_t}^{\,n-k}-\rho_{L_t}^{\,n-k}
=t^2h+o(t^2),
$$
and
$$
\rho_{K_t}^{\,n}-\rho_{L_t}^{\,n}
=\frac{n}{n-k}\,t^2\rho_K^k h+o(t^2),
$$
with both remainders uniform on \(\sph\).
	
We now apply the polar integration formulas~\eqref{eq:polar}.
Since every \(\SF\) has the same surface area, integrating the first
remainder over \(\SF\) gives \(o(t^2)\) uniformly in \(F\). We may
therefore take the maximum over \(F\) and obtain
\[
\max_{F\in G_{n,n-k}}
\bigl(|K_t\cap F|-|L_t\cap F|\bigr)
=\frac{t^2}{n-k}
\max_{F\in G_{n,n-k}}\int_{\SF}h\dd\sigma_F+o(t^2).
\]
Integrating the second expansion over \(\sph\) gives
\[
|K_t|-|L_t|
=\frac{t^2}{n-k}\int_{\sph}\rho_K^k h\dd\sigma+o(t^2).
\]
The radial functions of \(K_t\) and \(L_t\) converge uniformly to
\(\rho_K\), so both volumes tend to \(|K|>0\). Applying the mean
value theorem to \(s^{(n-k)/n}\) between these volumes yields
\[
|K_t|^{(n-k)/n}-|L_t|^{(n-k)/n}
=\frac{t^2}{n|K|^{k/n}}\int_{\sph}\rho_K^k h\dd\sigma+o(t^2).
\]
	
Choose \(v\in\sph\) with \(h(v)>0\). For any \(F\in G_{n,n-k}\)
containing \(v\), continuity and nonnegativity give
\(\int_{\SF}h\dd\sigma_F>0\). The leading coefficient in the maximum
section difference is therefore positive, and division gives
\[
\lim_{t\downarrow0}
\frac{|K_t|^{(n-k)/n}-|L_t|^{(n-k)/n}}
{\displaystyle\max_{F\in G_{n,n-k}}
\bigl(|K_t\cap F|-|L_t\cap F|\bigr)}
=
\frac{n-k}{n}\,
\frac{\displaystyle\int_{\sph}\rho_K^k h\dd\sigma}
{\displaystyle |K|^{k/n}
\max_{F\in G_{n,n-k}}\int_{\SF}h\dd\sigma_F}.
\]
By~\eqref{eq:volume-difference}, each quotient on the left is at most
\(\gamma_{n,k}^{\,k}\). Taking the supremum over \(K\) and \(h\) on the
right and combining with Step~1 gives
\begin{equation}\label{eq:angular-formula}
\gamma_{n,k}^{\,k}
=\frac{n-k}{n}\sup_{K,h}
\frac{\displaystyle\int_{\sph}\rho_K^k h\dd\sigma}
{\displaystyle |K|^{k/n}
\max_{F\in G_{n,n-k}}\int_{\SF}h\dd\sigma_F}.
\end{equation}
	
\textbf{Step~3.} For each fixed \(K\), we show that the supremum
over functions on the sphere in~\eqref{eq:angular-formula} equals the
supremum over densities on \(K\) in~\eqref{eq:density}.
	
Given a density \(f\) as in the proposition, define its weighted integral
along each ray by
\[
\begin{aligned}
h_f(\theta)&=\int_0^{\rho_K(\theta)} f(r\theta)r^{n-k-1}\dd r
=\rho_K(\theta)^{n-k}\int_0^1f(s\rho_K(\theta)\theta)s^{n-k-1}\dd s .
\end{aligned}
\]
The last integrand is continuous on \([0,1]\times\sph\), so \(h_f\) is continuous.
The evenness of \(f\) and \(\rho_K\) makes \(h_f\) even, and \(f\ge0\) gives \(h_f\ge0\).
Equations~\eqref{eq:polar-density} and~\eqref{eq:section-density} give
\[
\int_{K\cap F} f(x)\dd x=\int_{\SF}h_f\dd\sigma_F,\qquad
0<\int_K f(x)\dd x\le\int_{\sph}\rho_K^k h_f\dd\sigma.
\]
The inequality uses \(r^{n-1}=r^k r^{n-k-1}\) and \(r^k\le\rho_K(\theta)^k\).
In particular \(h_f\not\equiv0\). The density quotient in \eqref{eq:density}
has the same denominator and no larger numerator than its spherical counterpart.
	
For the converse fix \(K,h\) in \eqref{eq:angular-formula}. To make
\(r^k/\rho_K(\theta)^k\) close to one on each ray, place the mass in a thin
shell. For \(0<\delta<1\), choose a nonnegative continuous \(\psi_\delta\) supported
in \((1-\delta,1)\) and normalized by \(\int_0^1\psi_\delta(s)s^{n-k-1}\dd s=1\).
One may take a triangular bump on \([1-3\delta/4,1-\delta/4]\) and divide
by its positive weighted integral. Set
\[
f_\delta(r\theta)=\frac{h(\theta)}{\rho_K(\theta)^{n-k}}
\psi_\delta\!\Big(\frac r{\rho_K(\theta)}\Big),
\quad 0<r\le\rho_K(\theta),\qquad f_\delta(0)=0 .
\]
This is even and nonnegative. It is continuous away from the origin and vanishes
when \(r\le(1-\delta)\min\rho_K\), so it is continuous on all of \(K\).
The substitution \(s=r/\rho_K(\theta)\) gives, for either \(d=n-k\) or \(d=n\),
\[
\int_0^{\rho_K(\theta)} f_\delta(r\theta)r^{d-1}\dd r
=h(\theta)\rho_K(\theta)^{d-(n-k)}\int_0^1\psi_\delta(s)s^{d-1}\dd s .
\]
For \(d=n-k\) the last integral is one, so every section integral of \(f_\delta\)
equals the corresponding section integral of \(h\) on the sphere. For \(d=n\)
the last integral is between \((1-\delta)^k\) and \(1\), since it is the mean
of \(s^k\) against the probability measure \(\psi_\delta(s)s^{n-k-1}\dd s\).
Consequently
\[
\int_{K\cap F} f_\delta(x)\dd x=\int_{\SF}h\dd\sigma_F,\qquad
\int_K f_\delta(x)\dd x\ge(1-\delta)^k\int_{\sph}\rho_K^k h\dd\sigma>0.
\]
The strict inequality follows from positivity of \(h\) on a nonempty spherical open
set. For this fixed \(K,h\), the density supremum is therefore at least
\((1-\delta)^k\) times the spherical quotient. Letting \(\delta\downarrow0\)
and then taking the supremum over \(K,h\) proves \eqref{eq:density}.
\end{proof}

We now use \eqref{eq:density} to seek a lower bound. A ball illustrates
why the section integrals and the volume of the supporting body must be
controlled together. For the uniform probability density on \(r\B^n\),
the integral over every central section of codimension \(k\) is
\(q_r(0)=\omega_{n-k}/(\omega_n r^k)\). Hence
\[
|r\B^n|^{k/n}q_r(0)=\frac{\omega_{n-k}}{\omega_n^{(n-k)/n}},
\]
which is independent of \(r\). We will approximate the uniform measures
on balls of several radii by measures supported on finite sets of
points. At each radius we will control the number of points so that
their convex hull can be covered by sufficiently few balls.

The next lemma gives the approximation needed at a single radius. Since
convolution averages over translates, it must hold at every translation
and for every orientation of the section. Recall that \(q_r\) in
\eqref{eq:marginal} is the density of the projection of the uniform
probability measure \(\mu_r\) on \(r\B^n\).

\begin{lemma}\label{lem:ball-discretization}
Let \(n\ge3\), \(1\le k\le n/2\), \(r>0\), and \(0<\eta<1/2\).
Suppose that \(\varphi:\R^k\to[0,1]\) is supported in \(r\B^k\) and
satisfies \(\Lip(\varphi)\le L_0\sqrt n/r\), with \(L_0>0\) fixed.
There is an even discrete probability measure \(\nu\) supported on at most
\[
Cnk\,\eta^{-2}\log(Ck/\eta)
\]
points of \(r\B^n\) such that
\begin{equation}\label{eq:det-ball-discretization}
\Big|\int_{\R^n}\varphi(z-Ux)\dd\nu(x)
-\int_{\R^k}\varphi(z-y)q_r(y)\dd y\Big|\le\eta,
\qquad UU^\T=I_k,\ z\in\R^k.
\end{equation}
The constant \(C\) depends only on \(L_0\).
\end{lemma}

\begin{proof}
\textbf{Step~1.} We approximate the tests \(x\mapsto\varphi(z-Ux)\) by
averages of a finite family of functions on the ball.
For \(|z|>2r\) the test vanishes on \(r\B^n\), since \(|Ux|\le r\).
It suffices to treat \(|z|\le2r\). Write \(U=(u_{lj})\) and set
\[
 \delta=\frac{\eta}{8(1+L_0)\sqrt{nk}} .
\]
We use the coordinate rounding of Klartag and
Livshyts~\cite[Definition~2.2]{KLrounding} on each row of \(U\).
For a real number \(b=\delta(a+t)\) with \(a\in\mathbb Z\) and
\(0\le t<1\), give \(\delta a\) weight \(1-t\) and \(\delta(a+1)\)
weight \(t\). The two endpoints have weighted mean \(b\). Apply this
rule independently to the entries \(u_{lj}\), and denote the resulting
matrix by \(V\). Thus each row of \(V\) is a vertex of a lattice cube
containing the corresponding row of \(U\). Expectations in this step
are just finite averages over these vertices.

For any fixed \(x\in r\B^n\), the coordinate errors have mean zero
and absolute value at most \(\delta\). Independence makes the mixed
terms in each row vanish, so
\[
 \E |(V-U)x|^2
 =\sum_{l=1}^k\sum_{j=1}^n x_j^2\,\E(V_{lj}-u_{lj})^2
 \le k\delta^2 |x|^2 .
\]
Round each coordinate of \(z\) to the nearest multiple of \(r\delta\),
obtaining \(w\). Then \(|z-w|\le r\delta\sqrt k/2\). By the
Cauchy--Schwarz inequality and the Lipschitz bound, for every \(x\in r\B^n\),
\[
 \big|\varphi(z-Ux)-\E\varphi(w-Vx)\big|
 \le \frac{L_0\sqrt n}{r}
      \bigl(|z-w|+\sqrt{\E|(V-U)x|^2}\bigr)
 \le 2L_0\sqrt{nk}\,\delta\le\eta/4 .
\]
The weights depend only on \(U\), so this is a uniform approximation
on the ball. Each lattice matrix \(V\) serves only to define the test
\(x\mapsto\varphi(w-Vx)\).

\textbf{Step~2.} We count the possible rounded rows and translations by
comparing the volumes of cubes and balls.
Each row of \(U\) has length one. Each rounded row belongs to
\(\delta\mathbb Z^n\cap(1+\delta\sqrt n)\B^n\). Cubes of side
\(\delta\) centered at these lattice points have disjoint interiors
and lie in \(2\B^n\), since \(1+3\delta\sqrt n/2<2\).
The number of possible rows is consequently at most
\[
 \frac{2^n\omega_n}{\delta^n}
 \le \left(\frac{C}{\delta\sqrt n}\right)^n ,
\]
by~\eqref{eq:ball-volume}. Likewise the cubes of side \(r\delta\)
centered at the possible \(w\)'s lie in \(3r\B^k\), because
\(|z|\le2r\) and \(\delta\sqrt k<1\). Their number is at most
\((C/(\delta\sqrt k))^k\). Let \(H\) be the number of these choices of \((V,w)\). Combining the
\(k\) rows and the translation gives
\[
 H\le
 \Big(\frac{C}{\delta\sqrt n}\Big)^{nk}
 \Big(\frac{C}{\delta\sqrt k}\Big)^k,\qquad
 \log(2H)\le Cnk\log(Ck/\eta).
\]
Indeed, \(1/(\delta\sqrt n)=8(1+L_0)\sqrt k/\eta\) and
\(1/(\delta\sqrt k)=8(1+L_0)\sqrt n/\eta\). In the second factor,
\(k\log n\le nk\). The constant in the last bound depends only on \(L_0\).

\textbf{Step~3.} We choose points of the ball whose averages approximate
the integrals of these finitely many tests simultaneously.
Choose \(N\) points \(X_1,\ldots,X_N\) independently and uniformly in
\(r\B^n\). For a fixed test \(g(x)=\varphi(w-Vx)\), the centered values
\(Y_i=g(X_i)-\int g\dd\mu_r\) have mean zero and satisfy \(|Y_i|\le1\).
For \(0<t\le1\), the inequality \(e^s\le1+s+s^2\) for \(|s|\le1\)
gives \(\E e^{\pm tY_i}\le e^{t^2}\). Applying Markov's inequality to
\(e^{\pm t\sum_iY_i}\), multiplying the expectations by independence,
and taking \(t=\eta/4\), we get
\[
 \mathbb P\Big\{
 \Big|\frac1N\sum_{i=1}^N g(X_i)-\int g\dd\mu_r\Big|>\eta/2
 \Big\}
 \le 2e^{-N\eta^2/16}.
\]
For \(N=\lceil32\eta^{-2}\log(2H)\rceil\), the sum of these failure
probabilities is less than one. We can thus fix points \(x_1,\ldots,x_N\)
for which every rounded test has error at most \(\eta/2\).
Any weighted average of these tests has the same error bound. Returning
to the original test costs at most \(\eta/4\) in its integral over the
ball and \(\eta/4\) in its average over the selected points, by Step~1.
The total error is at most \(\eta\), for all \(U,z\).

Finally, give mass \(1/(2N)\) to each of \(x_i\) and \(-x_i\), counting
repetitions. The average against this even measure is the mean of the
two sample averages for \(U\) and \(-U\). Their integrals against
\(\mu_r\) agree by the symmetry of the ball, and
\(\int\varphi(z-Ux)\dd\mu_r(x)=\int\varphi(z-y)q_r(y)\dd y\)
by~\eqref{eq:marginal}. This proves~\eqref{eq:det-ball-discretization}
with at most \(2N\le Cnk\,\eta^{-2}\log(Ck/\eta)\) points.
\end{proof}

We apply Lemma~\ref{lem:ball-discretization} at successively smaller
radii. The bound on the number of atoms lets us estimate the volume of
the supporting convex body by Lemma~\ref{lem:convex-hull-covering}. The
error bound, which is uniform in the projection and the translation,
lets us estimate all section integrals of the resulting convolution.

\begin{lemma}\label{lem:small-sections}
There is an absolute constant \(C\) such that for every \(n\ge2\) and
\(1\le k\le n/2\) one can find an origin-symmetric convex body \(T\subset\R^n\)
and an even continuous probability density \(f\ge0\) on \(\R^n\) satisfying
\[
\supp f\subseteq T,\qquad |T|^{1/n}\le C,\qquad
\sup_{F\in G_{n,n-k}}\int_F f(x)\dd x\le\bigl(C\sqrt{k/n}\bigr)^k .
\]
\end{lemma}

\begin{proof}
For \(n=2\) we have \(k=1\). Take \(T=\B^2\) and
\(f(x)=(3/\pi)(1-|x|)_+\). Polar coordinates give
\(\int_{\R^2}f=6\int_0^1(1-t)t\dd t=1\). On every line \(F\) through
the origin, \(\int_F f=(6/\pi)\int_0^1(1-t)\dd t=3/\pi\), and
\(|T|^{1/2}=\sqrt\pi\). These are the required bounds with an absolute
constant. Henceforth \(n\ge3\), so that \(n-k\ge2\).

\textbf{Step~1.} We compute the effect of translating a ball on each
of its parallel sections.
For \(UU^\T=I_k\) and \(z\in\R^k\), write
\[
 F_z=\{x:Ux=z\}=U^\T z+\ker U .
\]
The section of \(a+s\B^n\) by \(F_z\) is a ball whose radius, when
\(|z-Ua|<s\), is \((s^2-|z-Ua|^2)^{1/2}\). Thus~\eqref{eq:marginal}
gives its normalized volume for every \(z\),
\[
 \frac{|(a+s\B^n)\cap F_z|}{\omega_n s^n}=q_s(z-Ua).
\]
If a density has integral \(Q(z)\) over \(F_z\), translating the
density by \(a\) changes that integral to \(Q(z-Ua)\).
Averaging the center \(a\) uniformly over \(s\B^n\) therefore gives
\[
 \int_{s\B^n}Q(z-Ua)\dd\mu_s(a)
   =\int_{\R^k}Q(z-y)q_s(y)\dd y .
\]

Choose the radii
\[
 R=\sqrt{n/k},\qquad J=\lceil\log_2 R\rceil,\qquad
 r_j=2^{-j}\sqrt n\,R, \quad 0\le j\le J.
\]
The largest ball satisfies \(q_{r_0}(0)\le R^{-k}\) by
Lemma~\ref{lem:ball-marginal}. The smallest has
\(\sqrt n/2<r_J\le\sqrt n\), so its volume to the power \(1/n\) is
bounded by~\eqref{eq:ball-volume}. We start with the smallest ball and
use the larger balls only to choose translation centers.

For the moment allow each translation center to range over its entire
ball. Start with \(Q_J=q_{r_J}\) and, in decreasing order of \(j\), set
\[
 Q_j(z)=\int_{\R^k}Q_{j+1}(z-y)q_{r_j}(y)\dd y,
 \qquad 0\le j<J.
\]
Thus \(Q_j(z)\) is the average normalized section volume of the
smallest ball after translations at radii \(r_{J-1},\ldots,r_j\).
Each \(Q_j\) has integral one. Since
\(Q_j=q_{r_j}*\cdots*q_{r_J}\), it is an average of translates
of \(q_{r_j}\). Averaging preserves upper bounds on both the supremum
norm and the Lipschitz constant. On taking convolutions, the support lies in the sum of the supports. Consequently, for \(j<J\),
\[
 \begin{aligned}
  0\le Q_{j+1}\le q_{r_{j+1}}(0), \qquad
  \supp Q_{j+1}\subseteq r_j\B^k, \qquad
  \Lip(Q_{j+1})\le
    \frac{2\sqrt{2n}}{r_j}\,q_{r_{j+1}}(0).
 \end{aligned}
\]
Here \(\sum_{i=j+1}^Jr_i\le r_j\) gives the support inclusion, and
Lemma~\ref{lem:ball-marginal} with \(r_{j+1}=r_j/2\) gives the last bound.

\textbf{Step~2.} We replace the continuous choices of translation
centers by finite sets and estimate the change in every section.
Apply Lemma~\ref{lem:ball-discretization} at radius \(r_j\) to
\(Q_{j+1}/q_{r_{j+1}}(0)\) with tolerance \(2^{-2k(j+1)}\).
Step~1 verifies its hypotheses with \(L_0=2\sqrt2\).
Let \(\nu_j\) be the resulting even probability measure and let \(N_j\)
be the number of its centers in \(r_j\B^n\). The point count in the
lemma gives \(N_j/n\le Ck\,2^{4k(j+1)}\log(Ck\,2^{2k(j+1)})\).
Taking logarithms and using \(\log(Ck)\le Ck(j+1)\), we obtain
\begin{equation}\label{eq:det-number-of-points}
 \log(1+N_j/n)\le Ck(j+1).
\end{equation}
For every \(U,z\), the estimate on the section averages is
\begin{equation}\label{eq:det-marginal-discrepancy}
 \begin{aligned}
 \Big|\int Q_{j+1}(z-Ua)\dd\nu_j(a)-Q_j(z)\Big|
 \le q_{r_{j+1}}(0)\,2^{-2k(j+1)}
 \le R^{-k}2^{-k(j+1)}.
 \end{aligned}
\end{equation}
The second inequality uses \(q_{r_{j+1}}(0)\le(2^{j+1}/R)^k\).

Let \(b=\mathbf{1}_{r_J\B^n}/(\omega_n r_J^n)\), the normalized
indicator of the smallest ball. Set \(g_J=b\) and, for \(j=J-1,\ldots,0\),
form the finite average
\[
 g_j(x)=\int g_{j+1}(x-a)\dd\nu_j(a).
\]
Each \(g_j\) is an even probability density and a finite weighted sum
of normalized indicators of translated copies of \(r_J\B^n\).
For a fixed \(U\), set \(P_j(z)=\int_{F_z}g_j\), with Lebesgue measure
on \(F_z\). The ball section formula gives \(P_J=Q_J\) at every \(z\).
Since \(F_z-a=F_{z-Ua}\), we have
\(P_j(z)=\int P_{j+1}(z-Ua)\dd\nu_j(a)\).
Subtracting \(Q_j(z)\), inserting the average of \(Q_{j+1}\), and
using~\eqref{eq:det-marginal-discrepancy} gives
\[
 \|P_j-Q_j\|_\infty
 \le \|P_{j+1}-Q_{j+1}\|_\infty+R^{-k}2^{-k(j+1)} .
\]
All the constants are independent of the orientation \(U\). Summing
from the smallest to the largest radius and using
\(\|Q_0\|_\infty\le q_{r_0}(0)\le R^{-k}\) yields
\[
 \sup_{UU^\T=I_k,\,z\in\R^k}\int_{F_z}g_0
 \le R^{-k}\Big(1+\sum_{j=0}^{J-1}2^{-k(j+1)}\Big)
 \le 2R^{-k}.
\]

\textbf{Step~3.} We obtain a continuous density from the translated
balls and cover their convex hull to bound its volume.
Average \(g_0\) over one further translation inside the smallest ball,
and set \(f(x)=\int g_0(x-a)b(a)\dd a\). Thus \(f\) is an even
probability density. To see its continuity, note that \(f\) is a
finite weighted sum of translates of
\[
 (b*b)(x)=
 \frac{|r_J\B^n\cap(x+r_J\B^n)|}{(\omega_n r_J^n)^2}.
\]
This intersection volume is continuous in \(x\). Indeed, as its center
moves to \(x\), the indicator of the translated ball converges at every
point off \(x+\partial(r_J\B^n)\). That sphere has volume zero, so
dominated convergence on \(r_J\B^n\) applies.
By Fubini's theorem, the section integrals of \(f\) are averages of
the parallel section integrals already estimated in Step~2. In particular,
\[
 p_U(z):=\int_{F_z} f=\int P_0(z-Ua)b(a)\dd a
 \le 2R^{-k}.
\]
The functions \(P_j\) are continuous in \(z\) by their recursion from
\(P_J=q_{r_J}\), so \(p_U\) is continuous as well. At \(z=0\) this
gives \(\int_{\ker U} f\le(2\sqrt{k/n})^k\).

Since the centers add when the balls are translated, we have
\[
 \supp f\subseteq 2r_J\B^n+\sum_{j=0}^{J-1}\supp\nu_j
 \subseteq T,\qquad
 T=2r_J\B^n+\sum_{j=0}^{J-1}\conv(\supp\nu_j).
\]
The finite sets of centers are symmetric about the origin, so \(T\)
is an origin-symmetric convex body. To cover \(T\), aim for radius
\(\sqrt n\,2^{-j/2}\) when covering the convex hull at level \(j\).
In Lemma~\ref{lem:convex-hull-covering} this is achieved by taking
\(\ell_j=\lceil n/(k2^j)\rceil\). Indeed, \(2^j<R\) for \(j<J\), so
\[
 \frac{n}{k2^j}\le\ell_j\le\frac{2n}{k2^j},\qquad
 \frac{r_j}{\sqrt{\ell_j}}\le\sqrt n\,2^{-j/2}.
\]
Thus \(\binom{N_j+\ell_j-1}{\ell_j}\) covering balls suffice at this
level. Using \(\ell_j!\ge(\ell_j/e)^{\ell_j}\),
\(1+N_j/\ell_j\le k2^j(1+N_j/n)\), and~\eqref{eq:det-number-of-points},
we get
\[
 \begin{aligned}
 \log\binom{N_j+\ell_j-1}{\ell_j}
 &\le \ell_j\,[1+\log(1+N_j/\ell_j)]\\
 &\le \frac{2n}{k2^j}
       [1+\log(k2^j)+\log(1+N_j/n)]
 \le Cn(j+1)2^{-j}.
 \end{aligned}
\]
Adding one covering ball at each level to \(2r_J\B^n\) gives a ball
of radius at most
\(\sqrt n(2+\sum_{j\ge0}2^{-j/2})<6\sqrt n\).
All such sums cover \(T\), and multiplying the numbers of choices
gives at most \(\exp(Cn\sum_{j\ge0}(j+1)2^{-j})=e^{4Cn}\) balls.
The sum of their volumes and~\eqref{eq:ball-volume} give
\[
 |T|^{1/n}\le 6e^{4C}\sqrt n\,\omega_n^{1/n}
 \le 6e^{4C}\sqrt{2\pi e}.
\]
This proves the volume bound for the same body that supports \(f\).
\end{proof}

The value \(p_U(0)\) also has a geometric interpretation. If \(X\) has
density \(f\) and \(F=\ker U\), then \(\operatorname{dist}(X,F)=|UX|\).
For every \(\varepsilon>0\),
\[
\frac{\mathbb P\{\operatorname{dist}(X,F)\le\varepsilon\}}
{\omega_k\varepsilon^k}
=\frac{1}{\omega_k\varepsilon^k}
\int_{|z|\le\varepsilon}p_U(z)\dd z .
\]
By continuity, this average tends to \(p_U(0)=\int_F f(y)\dd y\) as
\(\varepsilon\) decreases to zero.

We can now insert the density from Lemma~\ref{lem:small-sections} into
Proposition~\ref{prop:density}. This gives the lower bound for
\(k\le n/2\). For larger codimensions, \(\sqrt{n/k}\) is bounded and
it suffices to evaluate the volume-difference quotient for two
concentric Euclidean balls.

\begin{proof}[Proof of Theorem~\ref{thm:main}]
Suppose first that \(1\le k\le n/2\). Take \(T,f\) from
Lemma~\ref{lem:small-sections}, and choose \(C\ge1\) to dominate the
two constants in its estimates. Since \(f\) is a probability density
supported in \(T\), we have \(\int_T f=1\) and
\(\int_{T\cap F}f=\int_F f\) for every \(F\in G_{n,n-k}\).
Thus \eqref{eq:density} gives
\[
\begin{aligned}
\gamma_{n,k}^{\,k}
\ge \frac{n-k}{n}\,
\frac{1}{|T|^{k/n}
\displaystyle\max_{F\in G_{n,n-k}}\int_F f(x)\dd x}
\ge \frac{1}{2}\,
\frac{1}{C^k \bigl(C\sqrt{k/n}\bigr)^k}
=\frac12\Big(C^{-2}\sqrt{\frac nk}\Big)^k.
\end{aligned}
\]
Here \((n-k)/n\ge1/2\), and \(|T|^{1/n}\le C\) gives
\(|T|^{k/n}\le C^k\). Taking \(k\)-th roots and using
\(2^{-1/k}\ge1/2\) yields
\[
\gamma_{n,k}\ge (2C^2)^{-1}\sqrt{\frac nk}.
\]
	
For \(k=1\), the upper bound \(\gamma_{n,1}\le\sqrt n\) therefore gives
the sharp order of growth. The construction also yields a pair of bodies
violating any proposed absolute bound \(A>0\). Choose \(n\) so large
that \((2C^2)^{-1}\sqrt n>A\). For the pair \(T,f\) just used, the
density quotient in \eqref{eq:density}, including its prefactor, exceeds
\(A\). The spherical quotient for \(T,h_f\) is at least as large.
The proof of Proposition~\ref{prop:density} constructs bodies
\(L_t\subset K_t\) whose volume-difference quotients converge to this
spherical quotient. Consequently the quotient for \(L_t,K_t\) exceeds
\(A\) for all sufficiently small \(t>0\).
	
Now suppose that \(n/2<k<n\). Take \(K=\B^n\) and \(L=t\B^n\) with
\(0<t<1\). The left-hand side of \eqref{eq:volume-difference} is
\((1-t^{n-k})\omega_n^{(n-k)/n}\), and the difference on every
central section is \((1-t^{n-k})\omega_{n-k}\). Cancelling the
positive factor \(1-t^{n-k}\) gives
\[
\gamma_{n,k}^{\,k}\ge\frac{\omega_n^{(n-k)/n}}{\omega_{n-k}}.
\]
The same quotient is obtained in \eqref{eq:angular-formula} by taking
\(K=\B^n\) and \(h\equiv1\), since the two spherical integrals are then
\(n\omega_n\) and \((n-k)\omega_{n-k}\).
	
To bound this quotient, put \(a=2/\sqrt{2\pi e}\). By
\eqref{eq:ball-volume},
\[
\frac{\omega_n^{(n-k)/n}}{\omega_{n-k}}
=\Big(\frac{\omega_n^{1/n}}{\omega_{n-k}^{1/(n-k)}}\Big)^{n-k}
\ge \Big(a\sqrt{\frac{n-k}{n}}\Big)^{n-k}.
\]
Since \(0<a<1\) and \(0<(n-k)/k<1\), we have \(a^{(n-k)/k}\ge a\).
For the remaining factor after taking \(k\)-th roots, the elementary
bound \(\log(1+s)\le s\) for \(s\ge0\) gives
\[
\frac{n-k}{2k}\log\frac{n}{n-k}
=\frac{n-k}{2k}\log\Big(1+\frac{k}{n-k}\Big)\le\frac12.
\]
It follows that \(((n-k)/n)^{(n-k)/(2k)}\ge e^{-1/2}\), and hence
\[
\gamma_{n,k}\ge a e^{-1/2}
\ge \frac{a}{\sqrt{2e}}\sqrt{\frac nk},
\]
where the last inequality uses \(n/k<2\).
Taking \(c=\min\{(2C^2)^{-1},a/\sqrt{2e}\}\) proves
\eqref{eq:bounds} for every \(1\le k<n\).
\end{proof}

\noindent\textbf{Declaration on the use of AI.}
OpenAI was used to assist with manuscript revision and the development and verification of the numerical code. The author retains full responsibility for the mathematical arguments, numerical results, and final content of the manuscript.

\end{document}